\documentstyle[12pt]{article}
\begin{document}
% \large
\begin{center}

{\bf EXACT EXPONENTIAL BOUNDS FOR THE RANDOM }

\vspace{3mm}

  {\bf FIELD MAXIMUM DISTRIBUTION VIA THE }\\

\vspace{3mm}

  {\bf MAJORING MEASURES (GENERIC CHAINING)}\\

\vspace{3mm}

{\sc By Ostrovsky E., Rogover E.}\\

\vspace{3mm}
{\it Department of Mathematics and Statistics, Bar-Ilan University, 
59200, Ramat Gan, Israel.}\\
e-mail: \ galo@list.ru \\

{\it Department of Mathematics and Statistics, Bar-Ilan University, 
59200, Ramat Gan, Israel.}\\
e - mail: \ rogovee@gmail.com \\

\end{center}
\vspace{3mm}

  In this paper non-asymptotic exact exponential estimates
 are derived for the tail of maximum's distribution
of random field in the terms of majoring measures or, equally,
 generic chaining.\\

 {\it Key words:} Majoring measures, generic chaining, random fields,
 exponential estimations, entropy,
 moment, Banach spaces of random variables, tail of distribution.\\

{\it Mathematics Subject Classification (2000):} primary 60G17; \ secondary
 60E07; 60G70.\\

\vspace{3mm}

{\bf 1. Introduction. Notations. Statement of problem.}
 Let $ (\Omega,F,{\bf P} ) $ be a probability space, $ \Omega = \{\omega\}, \
T = \{t \} $ be arbitrary set, $ \xi(t), \ t \in T $ be centered:
$ {\bf E} \xi(t) = 0 $ separable random field (or process). For arbitrary
subset $ S \subset T $ we denote
$$
Q(S,u) = {\bf P}( \sup_{t \in S} \xi(t) > u), \ u \ge 2. \eqno(1)
$$
$$
Q_+(S,u) = {\bf P}( \sup_{t \in S} |\xi(t)| > u), \ u \ge 2. \eqno(2)
$$

{\bf Our aim is obtaining an exponentially exact as $ u \to \infty $
 estimation for the probability $ Q(u) \stackrel{def}{=} Q(T,u) $ in the
 terms of "majoring measures" or equally  in the
 terms of "generic chaining". } \par

  Definitions and some important results about $ {\bf E} \sup_{t \in T}
\xi(t) $ in the terms of majoring measures see, for example, in [2], [3],
[6], p. 309 - 330, [10], [11], [12], [13]. In the so-called "entropy"
terms this problem was considered in [3], [4]. See also [8]. \par
 Note that the "majoring measures" method is more general in comparison
to the entropy technique ([6], p. 309 - 330, [10], [12]). \par

\vspace{3mm}

{\bf 2. Auxiliary facts.}
 In order to formulate our result, we need to introduce some addition
notations and conditions. Let $ \phi = \phi(\lambda), \lambda \in (-\lambda_0, \lambda_0), \ \lambda_0 = const \in (0, \infty] $ be some even strong convex which takes positive values for positive arguments twice continuous
differentiable function, such that
$$
 \phi(0) = 0, \ \phi^{//}(0) > 0, \ \lim_{\lambda \to \lambda_0} \phi(\lambda)/\lambda = \infty. \eqno(3)
$$
 We denote the set of all these function as $ \Phi; \ \Phi =
\{ \phi(\cdot) \}. $ \par
 We say that the {\it centered} random variable (r.v) $ \xi = \xi(\omega) $
belongs to the space $ B(\phi), $ if there exists some non-negative constant
$ \tau \ge 0 $ such that

$$
\forall \lambda \in (-\lambda_0, \lambda_0) \ \Rightarrow
{\bf E} \exp(\lambda \xi) \le \exp[ \phi(\lambda \ \tau) ]. \eqno(4).
$$
 The minimal value $ \tau $ satisfying (4) is called a $ B(\phi) \ $ norm
of the variable $ \xi, $ write
 $$
 ||\xi||B(\phi) = \inf \{ \tau, \ \tau > 0: \ \forall \lambda \ \Rightarrow
 {\bf E}\exp(\lambda \xi) \le \exp(\phi(\lambda \ \tau)) \}.
 $$
 This spaces are very convenient for the investigation of the r.v. having a
exponential decreasing
 tail of distribution, for instance, for investigation of the limit theorem,
the exponential bounds of distribution for sums of random variables,
non-asymptotical properties, problem of continuous of random fields,
study of Central Limit Theorem in the Banach space etc.\par

  The space $ B(\phi) $ with respect to the norm $ || \cdot ||B(\phi) $ and
ordinary operations is a Banach space which is isomorphic to the subspace
consisted on all the centered variables of Orlich's space $ (\Omega,F,{\bf P}), N(\cdot) $ with $ N \ - $ function

$$
N(u) = \exp(\phi^*(u)) - 1, \ \phi^*(u) = \sup_{\lambda} (\lambda u -
\phi(\lambda)).
$$
 The transform $ \phi \to \phi^* $ is called Young-Fenchel transform. The proof of considered assertion used the properties of saddle-point method and theorem
 of Fenchel-Moraux:
$$
\phi^{**} = \phi.
$$

 The next facts about the $ B(\phi) $ spaces are proved in [4], [8, p. 19-40]:

$$
{\bf 1.} \ \xi \in B(\phi) \Leftrightarrow {\bf E } \xi = 0, \ {\bf and} \ \exists C = const > 0,
$$

$$
U(\xi,x) \le \exp(-\phi^*(Cx)), x \ge 0,
$$
where $ U(\xi,x)$ denotes in this article the {\it tail} of
distribution of the r.v. $ \xi: $

$$
U(\xi,x) = \max \left( {\bf P}(\xi > x), \ {\bf P}(\xi < - x) \right),
\ x \ge 0,
$$
and this estimation is in general case asymptotically exact. \par
 Here and further $ C, C_j, C(i) $ will denote the non-essentially positive
finite "constructive" constants.\par
 More exactly, if $ \lambda_0 = \infty, $ then the following implication holds:

$$
\lim_{\lambda \to \infty} \phi^{-1}(\log {\bf E} \exp(\lambda \xi))/\lambda =
K \in (0, \infty)
$$
if and only if

$$
\lim_{x \to \infty} (\phi^*)^{-1}( |\log U(\xi,x)| )/x = 1/K.
$$
 Here and further $ f^{-1}(\cdot) $ denotes the inverse function to the
function $ f $ on the left-side half-line $ (C, \infty). $ \par

The function $ \phi(\cdot) $ may be "constructive" introduced by the formula
$$
\phi(\lambda) = \phi_0(\lambda) \stackrel{def}{=} \log \sup_{t \in T}
 {\bf E} \exp(\lambda \xi(t)), \eqno(5)
$$
 if obviously the family of the centered r.v. $ \{ \xi(t), \ t \in T \} $ satisfies the {\it uniform } Kramer's condition:
$$
\exists \mu \in (0, \infty), \ \sup_{t \in T} U(\xi(t), \ x) \le \exp(-\mu \ x),
\ x \ge 0.
$$
 In this case, i.e. in the case the choice the function $ \phi(\cdot) $ by the
formula (5), we will call the function $ \phi(\lambda) = \phi_0(\lambda) $
a {\it natural } function. \par
 {\bf 2.} We define $ \psi(p) = p/\phi^{-1}(p), \ p \ge 2. $
 Let us introduce a new norm (the so-called "moment norm")
on the set of r.v. defined in our probability space by the following way: the space $ G(\psi) $ consist, by definition, on all the centered r.v. with finite norm

$$
||\xi||G(\psi) \stackrel{def}{=} \sup_{p \ge 2} |\xi|_p/\psi(p), \ |\xi|_p =
{\bf E}^{1/p} |\xi|^p. \eqno(6)
$$

 It is proved that the spaces $ B(\phi) $ and $ G(\psi) $ coincides:$ B(\phi) =
G(\psi) $ (set equality) and both
the norm $ ||\cdot||B(\phi) $ and $ ||\cdot|| $ are equivalent: $ \exists C_1 =
C_1(\phi), C_2 = C_2(\phi) = const \in (0,\infty), \ \forall \xi \in B(\phi) $

$$
||\xi||G(\psi) \le C_1 \ ||\xi||B(\phi) \le C_2 \ ||\xi||G(\psi).
$$

{\bf 3.} The definition (6) is correct still for the non-centered random
variables $ \xi.$ If for some non-zero r.v. $ \xi \ $ we have $ ||\xi||G(\psi) < \infty, $ then for all positive values $ u $

$$
{\bf P}(|\xi| > u) \le 2 \ \exp \left( - u/(C_3 \ ||\xi||G(\psi)) \right).
\eqno(7)
$$
and conversely if a r.v. $ \xi $ satisfies (7), then $ ||\xi||G(\psi) <
\infty. $ \par
{\sc We suppose in this article that there exists a function $ \phi \in \Phi $
such that} $ \forall t \in T \ \Rightarrow \ \xi(t) \in B(\phi) $ and
$$
  \sup_t [ \ ||\xi(t) \ ||B(\phi)] < \infty,
$$
or equally $ {\bf E} \xi(t) = 0, \ t \in T, $ and for all non-negative values
$ x $

$$
\sup_t \max \left[ {\bf P}\left( \xi(t) > x \right), \ {\bf P}
\left( \xi(t) < - x \right) \right] \le \exp \left(-\phi^*(x) \right). \eqno(8)
$$
 Note that if for some $ C = const \in (0,\infty) $
$$
Q_+(T,u) \le \exp \left(-\phi^*(C u) \right),
$$
 then the condition (8) is satisfied (the "necessity" of the condition (8)). \par
 M.Talagrand [10] – [13], W.Bednorz [2], X. Fernique [3] etc. write
instead our function $ \exp \left(-\phi^*(x) \right) $ some Young's
function $ \Psi(x) $ and used as a rule a function $ \Psi(x) = \exp(-x^2/2) $
( the so-called "subgaussian case").\par

 Without loss of generality we can and will suppose

$$
\sup_{t \in T} [ \ ||\xi(t) \ ||B(\phi)] = 1,
$$
(this condition is satisfied automatically in the case of natural choosing
of the function $ \phi: \ \phi(\lambda) = \phi_0(\lambda) \ ) $
and that the metric space $ (T,d) $ relatively the so-called
{\it natural} distance (more exactly, semi-distance)

$$
d(t,s) \stackrel{def}{=} ||\xi(t) - \xi(s)|| B(\phi)
$$
is complete. \par
 Recall that the semi-distance $ d = d(t,s), \ s,t \in T $ is, by definition,
non-negative symmetrical numerical function, $ d(t,t) = 0, \ t \in T, $
satisfying the triangle inequality, but the equality $ d(t,s) = 0 $
does not means (in general case) that $ s = t. $ \par
 For example, if $ \xi(t) $ is a centered Gaussian field with covariation function
 $ D(t,s) = {\bf E} \xi(t) \ \xi(s), $ then
$ \phi_0(\lambda) = 0.5 \ \lambda^2, \ \lambda \in R, $ and $ d(t,s) = $

$$
||\xi(t) - \xi(s)||B(\phi_0) = \sqrt{ \bf{Var} [ \xi(t) - \xi(s) ]} =
\sqrt{ D(t,t) - 2 D(t,s) + D(s,s) }.
$$

 There are many examples of {\it martingales}, e.g., in the article [7],
 $ ( \xi(n), F(n) ), \ T = \{1,2,3, \ldots, n, \ldots \} $
satisfying the following modification (8a) of the condition (8):

$$
\sup_n U(\xi(n)/\sigma(n), \ x) \le \exp \left[-\phi^*(x) \right], \eqno(8a),
$$
in particular, there are many examples with

$$
\phi^*(x) \sim x^r \ L^{1/r}(x)/r, \ x \to \infty; \ r = const \ge 1,
\eqno(9)
$$
where as usually $ f(x) \sim g(x), \ x \to \infty $ denotes
$$
\lim_{x \to \infty} f(x)/g(x) = 1;
$$
and with
$$
n^{\beta} \ L_1(n) \le \sigma(n) \le n^{\beta} \ L_2(n), \ \beta = const > 0,\eqno(10)
$$

$ L_1(x), L_2(x), L(x) $ are some positive continuous {\it slowly varying} as
$ x \to \infty $ functions,

$$
\sigma(n) = \sqrt{ {\bf Var}( \xi(n) ) }.
$$

 It is known ( [4], [8], p. 22 - 25) that (9) is equivalent in the case
$ r > 1 \ $ (under some simple assumption) to the following equality:

$$
\lambda \to \infty \ \Rightarrow
\phi(\lambda) \sim \lambda^s \ L^{- 1/s} (\lambda^s)/s, \ s = r/(r-1).
$$

 Let us introduce for any subset $ V, \ V \subset T $ the so-called
{\it entropy } $ H(V, d, \epsilon) = H(V, \epsilon) $ as a logarithm
of a minimal quantity $ N(V,d, \epsilon) = N(V,\epsilon) = N $
of a balls $ S(V, t, \epsilon), \ t \in V: $
$$
S(V, t, \epsilon) \stackrel{def}{=} \{s, s \in V, \ d(s,t) \le \epsilon \},
$$
which cover the set $ V: $
$$
N = \min \{M: \exists \{t_i \}, i = 1,2,…, M, \ t_i \in V, \ V
\subset \cup_{i=1}^M S(V, t_i, \epsilon ) \},
$$
and we denote also
$$
H(V,d,\epsilon) = \log N; \ S(t_0,\epsilon) \stackrel{def}{=}
 S(T, t_0, \epsilon), \ H(d, \epsilon) \stackrel{def}{=} H(T,d,\epsilon).
$$
 It follows from Hausdorf's theorem that
$ \forall \epsilon > 0 \ \Rightarrow H(V,d,\epsilon)< \infty $ iff the
metric space $ (V, d) $ is precompact set, i.e. is the bounded set with
compact closure.\par

 Now we recall, modify and rewrite some definition of "generic chaining"
theory, belonging to X.Fernique [3] and M.Talagrand [10] - [13]. Let the ball
$ S(T,t_0, \delta) = S(t_0, \delta), \ t_0 \in T, \ \delta \in (0,1] $
be a given. The {\it sequence} $ R $ of finite subsets of $ S(t_0, \delta) $
$ T_m, m = 0,1,2, \ldots, T_m \subset S(t_0, \delta), R = \{ T_m \} $
such that $ T_0 = \{ t_0 \}; $ here
and further a symbol $ | V| $ will denote the number of elements of
a finite set $ V: \ |V| = card(V), $
and such that the set $ \cup_{n=0}^{\infty} T_n $ is dense in $ T $
with respect to the semi-distance $ d, $ is called {\it generic chaining } of
$ S(t_0,\delta). $ The set of all generic chaining we will denote $ W: $
$$
W = W(S(t_0, \delta)) = W(t_0, \delta) \stackrel{def}{=} \{ R \}.
$$
 For any element $ t \in S(t_0, \delta) $ we denote arbitrary, but
fixed (non-random) element $ \pi_n(t)$ of a subset of
$ T_n: \pi_n(t) \in T_n $ such that
$$
d(t, \pi_n(t)) = \min_{s \in T_n} d(t,s). \eqno(11)
$$
 Let $ \gamma = \{ \gamma_n \}, \ n = 1,2, \ldots $ be arbitrary fixed non-random sequence of a positive numbers such that
$$
\sum_n 1/\gamma_n = 1, \ \gamma_1 \ge 3; \eqno(12)
$$
for example, $ 1/\gamma_n = \rho^{n-1} (1 - \rho), \ \rho = const \in (2/3,1). $
Let us introduce the following important function
$$
 L(t_0, \delta, R, \gamma) = L(t_0, \delta, R) = \sup_{t \in B(t_0,\delta) } \sum_{m=1}^{\infty} d(\pi_m(t), \pi_{m-1}(t))/\gamma_m. \eqno(13)
$$
 We will consider only the so-called "admissible" random fields (in the terms of
M.Talagran) $ \xi(\cdot), $
i.e. which satisfied the following conditions. \par
  Let us denote
$$
K(\xi, \phi, \delta) = K(\delta) = \inf_{R \in A} \ \inf_{ \gamma} \ 
\sup_{t_0 \in T} L(t_0, \delta, R, \gamma),
$$
 if the set $ A $ is not empty, and $ K(\delta) = + \infty $ in the other
case. \par
 {\it The following conclusions will be interest only in the case if for
some function } $ \phi(\cdot) \in \Phi, $ {\it for example for the natural function} $ \phi_0(\cdot) $

$$
\lim_{\delta \to 0+} K(\xi, \phi, \delta) = 0. \eqno(14)
$$
 {\sc We will suppose moreover that the condition (14), which will
called the "uniform generic chaining condition", write: $ \xi(\cdot), \phi
\in UA, $ is satisfied.}\par

 Let us introduce also the events $ D, \ E(n), n \ge 1 $ as follow:
$ E(n) = E(u; n, t_0, \delta, R) = $

$$
\cup_{t \in S(t_0, \delta)} [ \xi(\pi_n(t)) - \xi(\pi_{n-1}(t)) >
u \ d(\pi_n(t), \pi_{n-1}(t))/\gamma_n ] =
$$

$$
\cup_{t \in T_{n - 1} } [ \xi(\pi_n(t)) - \xi(\pi_{n -1 }(t)) >
u \ d(\pi_n(t), \pi_{n - 1}(t))/\gamma_n] =
$$

$$
 \{ \omega: \max_{t \in T_{n - 1} }
 \frac{ \xi(\pi_{n - 1} (t) - \xi(\pi_n(t))}{d(\pi_{n - 1}(t), \pi_n(t)) }
> \frac{u}{\gamma_n} \},
$$
if we define $ 0/0 = 0 $ (in the case if $ d(\pi_{n-1}(t),\pi_n(t)) = 0); $

$$
 D = D(u; t_0, \delta, R) = \cup _{n=1}^{\infty} E(u; n, t_0, \gamma, R).
$$

 We denote also

$$
Z_n = Z_n(u; t_0, \delta, \gamma, R) = {\bf P} [ E(u; n, t_0, \gamma, R)],
$$

$$
Y(u) = Y(u; t_0, \delta,R) = {\bf P} [ D(u; t_0, \delta, R ) ].
$$

 It is evident that

$$
 Z_n \le |T_n| \ |T_{n-1}| \exp \left(-\phi^*(u/\gamma_n) \right),
$$

$$
Y(u) \le \sum_{n=1}^{\infty} |T_n| \ |T_{n-1}| \exp \left(-\phi^*(u/\gamma_n) \right) \stackrel{def}{=} X(u; t_0, \delta, \gamma, R), \eqno(15)
$$

since
$$
{\bf P} [ \xi(\pi_n(t))- \xi(\pi_{n-1}(t)) > u \ d( \pi_n(t), \pi_{n-1}(t) ) ] \le
\exp \left( - \phi^*(u) \right), \ u > 0.
$$

 {\it The random field $ \xi(t) $ and the function $ \phi(\cdot) \ $ satisfies,
by definition, the uniform generic chaining condition, and write $ \xi(\cdot) \in A, $ or more simple: there exists the {\it set} of generic chaining $ R $ (depending on the $ \xi(\cdot) ) $ belonging to $ A, \ R \in A, $ if for
all $ \delta \in (0,1) $ and for arbitrary ball $ S(t_0, \delta) $ there exists
(for some sequence $ \{\gamma \} ) $
a generic chaining $ R $ in $ S(t_0, \delta) $ for which }

$$
Y(u; t_0, \delta, \gamma, R) \le \exp \left(-\phi^*(u/2) \right),
\eqno(16)
$$

{\it if, for example,}

$$
X(u; t_0, \delta, \gamma, R) \le \exp \left(-\phi^*(u/2) \right).\eqno(16a)
$$

 The existence of such a generic chaining it follows from our next assumptions.\par

{\bf Lemma 1.} {\it We have under the conditions (14) and (16) (or (16a))
for all the values} $ \delta \in (0,1]: $

$$
\sup_{t_0 \in T} {\bf P}
\left[ \sup_{t \in S(t_0, \delta)} (\xi(t) - \xi(t_0)) > u K(\delta) \right]
\le \exp \left(-\phi^*(u/2) \right). 
$$

{\bf Proof.} The proof of this assertion is alike to the original proof
of Talagran ([10], [12, chapter 1, pp. 9 - 14]) for the probability $ Q(u) = Q(T; u) $ estimation. Namely, let $ t_0 $ be arbitrary element of $ T, \
\delta \in (0,1]. $ Let also
$ R = \{ T_0, T_1, T_2, \ldots \}, \ R \in W $ be arbitrary chaining
into the ball $ S(t_0, \delta). $ We rewrite the Talagran's decomposition
([12], chapter 1, p. 10) for the ball $ S(t_0, \delta): $

$$
\xi(t) - \xi(t_0) = \sum_{n=1}^{\infty} [ \xi(\pi_n(t)) -
\xi(\pi_{n-1}(t) ) ].
$$
 Recall that $ \pi_0(t) = t_0 $ and that $ \forall t \in T $
$$
\lim_{n \to \infty} \pi_n(t) = t, \ \lim_{n \to \infty} \xi(\pi_n(t)) =
\xi(t)
$$
 in the sense of convergence in probability. \par

 We get analogously to the works [10], [11] and taking into account the
inclusion $ R \in A: $

$$
G(u; t_0, \delta) \stackrel{def}{=}{\bf P} \left( \sup_{t \in S(t_0,\delta)}
(\xi(t) - \xi(t_0)) > u \ K(\delta) \right) \le
$$
$$
Y(u; t_0, \delta, \gamma, R) \le \exp \left(-\phi^*(u/2) \right). \eqno(17)
$$

%$$
%\sum_{m=1}^{\infty} |T_m| \ |T_{m-1} | \exp(- \phi^*(u/\gamma_m) ) \le
%\exp(-\phi^*(u/2)). \eqno(17)
%$$

  Note that it follows from conclusion of Lemma 1 the {\it continuity} of
$ \xi(t) $ with probability one in the semi-distance $ d: $
$$
{\bf P}( \xi(\cdot) \in C(T,d) ) = 1;
$$
$ C(T,d) $ denotes as usually the space of all continuous with
respect to the semi-distance $ d $ functions $ f: T \to R. $ \par

 The conditions (14) and (16) is equivalent to the so-called condition of
the "uniform
 convergence of the majoring integral", see [10], [11]. \par

\vspace{3mm}

{\bf 3. Main result.} Let us denote for $ h \in \left(0, \sup_{\delta \in
(0,1 )} K(\delta) \right) \stackrel{def}{=} (0,K_0) $

$$
K^{-1}(h) = \inf \{\delta, \ \delta \in (0,1), \ K(\delta) \ge h \},
$$

$$
\Delta(C,u) = \Delta_{\phi}(u) = K^{-1}\left[ 0.5 \ C/ \left(u \ d
\phi^*(u)/d u) \right) \right],
$$
where $ d/du $ denotes the {\it right } derivative; it is obvious that the
derivative $ d \phi^*(u)/d u $ there exists, is continuous and
the function $ u \to \Delta(u) $
tends monotonically to zero as $ u \to \infty.$ Therefore, for arbitrary constant
$ C \in [1, \infty) $ there is a positive value $ u_0 = u_0(C), $
for which $ u \ge u_0 \ \Rightarrow \Delta(C, u) \le 0.5 K_0. $

{\bf Theorem 1.} {\it Suppose for any function } $ \phi(\cdot) \in \Phi $
$$
\lim_{\delta \to 0+} K(\xi(\cdot), \phi, \delta) = 0
$$
{\it and suppose the condition (16), or, more generally, (16a) is also
 satisfied.} \par
{\it Then for arbitrary constant } $ C \in (0,\infty) $
{\it and for all the values } $ u \ge u_0(C) $

$$
Q(u) \le [\exp(C) + 1) ] \ N(T,d, C \Delta_{\phi}(u)) \ \exp(-\phi^*(u)).
\eqno(18)
$$

 As a consequence:	
$$
Q_+(u) \le 2 \ [\exp(C) + 1) ] \ N(T,d, C \Delta_{\phi}(u)) \
\exp(-\phi^*(u)).
$$

{\bf Proof.} {\sc Step 1.} Let $ C $ be arbitrary positive constant,

$$
 u \ge u_0(C), \ \delta_0 = \delta_0(u) \stackrel{def}{= }
K^{-1}(0.5 \ C \ \Delta(u)).
$$

 We consider at first the probability

$$
Q(W,u) = {\bf P}\left(\sup_{t \in S(t_0,\delta_0)} \xi(t) > u \right),
$$
$ W = S(t_0,\delta_0). $ Denote $ \beta = C \ \Delta(u), \ \alpha = 1 - \beta;$
then $ \alpha, \beta > 0, \ \alpha + \beta = 1. \ $ We obtain:

$$
Q(W,u)\le {\bf P}(\xi(t_0) > \alpha \ u) + {\bf P}(\sup_{t \in W}
(\xi(t) - \xi(t_0)) > \beta \ u) \stackrel{def}{=} \eqno(19)
$$
$ I_1 + I_2. $ For the first member is true the simple estimation:
$$
I_1 \le \exp \left( - \phi^*(\alpha \ u ) \right).
$$
 As long as the function $ x \to \phi^*(x) $ is convex and twice
differentiable,

$$
\phi^*(\alpha \ u) \ge \phi^*(u) - (\phi^*)^/(u) \ C \ u \ \Delta(u) =
\phi^*(u) - C;
$$
 therefore

$$
I_1 \le \exp(C) \ \exp(-\phi^*(u)).
$$
 Further, since

$$
K(\delta_0) \le 0.5 \ C \ \Delta(u),
$$
we conclude using the inequality (17)

$$
I_2 \le \exp(-\phi^*(u)).
$$
 Summing, we receive:

$$
Q(W,u) \le C_1 \exp \left(-\phi^*(u) \right), \ C_1 = 1 + \exp C.
$$

\vspace{2mm}

{\sc Step 2.} Let $ \epsilon = C \ \Delta(u) $ and $ \{ t_i \}, \
i = 1,2,\ldots, N, $ where $ N = N(T,d,\epsilon) $ be a centers of
a balls $ B(T,d,\epsilon) $ forming a {\it minimal} (not necessary to be
unique) $ \epsilon \ - $ net of $ T $ with respect to the semi-distance
$ d. $\par
 Since the probability $ Q(S,u) $ has a property
$$
Q(S_1 \cup S_2,u) \le Q(S_1,u) + Q(S_2,u), \ S_1, S_2 \subset T,
$$

we conclude:

$$
Q(T,u) \le \sum_{i=1}^N Q(S(t_i, C \Delta(u)) ).
$$
 The last probabilities was estimated in (18).\par
 The low bounds for probabilities $ Q(T,u) , \ Q_+(T,u) $ was obtained
in ([8], 105 - 117); see also [9].\par
{\bf Corollary.} We explain here the exponential exactness of the estimation
of theorem 1. \par
  In many practical cases (statistics, method Monte-Carlo etc.) the entropy
  $ N(T, d, C \Delta_{\phi}(u) $ satisfies the inequality: $ \forall
  \epsilon \in (0,1/2) \ \exists U = U(\epsilon) \in (0, \infty) \ \Rightarrow
  \forall u \ge U(\epsilon) $

  $$
  N(T, d, C \Delta_{\phi}(u) \le \exp \left( \phi^*( \epsilon \ u ) \right),
  $$
   for example,

  $$
  N(T, d, C \Delta_{\phi}(u)) \le C (u + 1)^{\kappa}, \ \kappa \in (0, \infty),
   u \ge 0.
  $$

   Therefore in this cases

 $$
 Q(T,u) \le C_1(\epsilon) \exp \left( - \phi^*(( 1 - \epsilon))u \right).
 $$
  But there exists a {\it random variable }
 $ \xi, \ \xi \in B(\phi), \ ||\xi||B(\phi) = 1 $
  for which for $ u \ge U(\epsilon) $

 $$
 {\bf P}(\xi > u) \ge C_2(\epsilon) \exp \left( - \phi^*(( 1 + \epsilon))u \right).
 $$

\vspace{3mm}

{\bf 4. Exponential bounds for the sums of random fields.}
 Let $ \{ \xi_i(t) \}, \ i = 1,2,\ldots $ be an independent copies of
$ \xi(t), $
$$
 \eta_n(t) = n^{-1/2} \sum_{i=1}^n \xi_i(t),
$$

$$
Q_n(S,u) = {\bf P} \left(\sup_{t \in S} \eta_n(t) > u \right), \ Q_n(u)
= Q_n(T,u),
$$

$$
Q_{\infty} (S,u) = \sup_n Q_n(S,u), \ Q_{\infty}(u) = Q_{\infty}(T,u).
$$

{\bf We obtain in this section using (18) the
exponentially exact as $ u \to \infty $ in the aforementioned sense
uniform and non-uniform estimations for the probabilities
$ Q_n(u), Q_{\infty}(T,u) $ again in the terms of "generic chaining". } \par
 In the "entropy" terms this estimations are obtained in [1], [5].\par
 Let us denote for $ \lambda \in (-\lambda_0, \lambda_0) $
$$
\phi_n(\lambda) = n \ \phi(\lambda/\sqrt{n}), \ \zeta(\lambda) = \sup_n
\phi_n(\lambda),
$$
and introduce some new semi-distances:
$$
d_n(t,s) = ||\xi(t) - \xi(s)||B(\phi_n),
$$
$$
r(t,s) = ||\xi(t) - \xi(s)||B(\zeta).
$$
 As long as there exists a limit $ \lim_{n \to \infty} n \
\phi(\lambda/\sqrt{n}) = \sigma^2 \lambda^2/2, \ \sigma^2 = const
\in (0,\infty), $ we conclude that the function $ \zeta(\cdot) $
exists, is non-trivial and convex. \par

{\bf Theorem 2.} \par

{\bf A.} {\it Suppose for some function } $ \phi(\cdot) \in \Phi $

$$
\lim_{\delta \to 0+} K(\eta_n(\cdot), \phi_n, \delta) = 0.
$$

{\it Then for arbitrary constant } $ C \in (0,\infty) $
{\it and for all the values } $ u \ge u_0(C) $ {\it the following inequality
holds:}
$$
Q_n(u) \le [\exp(C) + 1) ] \ N(T,d_n, C \Delta_{\phi_n}(u)) \
 \exp(-\phi_n^*(u)).
$$

 {\bf B.} {\it Suppose for some function } $ \phi(\cdot) \in \Phi $

$$
\lim_{\delta \to 0+} \sup_n K(\eta_n(\cdot), \zeta, \delta) = 0.
$$

{\it Then for arbitrary constant } $ C \in (0,\infty) $
{\it and for all the values } $ u \ge u_0(C) $

$$
Q_{\infty}(u) \le [\exp(C) + 1) ] \ N(T, r, C \Delta_{\zeta}(u)) \
\exp(-\zeta^*(u)).
$$

  The conclusion of theorem 2 it follows trivially from the theorem 1 and
the following elementary fact: if $ \theta \in B(\phi), \ \phi \in \Phi, $
and $ \theta(i) $ are independent copies of $ \theta, $
$$
 \nu_n \stackrel{def}{=} n^{-1/2}\sum_{i=1}^n \theta(i),
$$
 then

$$
{\bf E} \exp(\lambda \ \nu_n) \le \exp(\phi_n(\lambda)).
$$

 Note that under the conditions of theorem 2 {\bf B} the sequence of the
random fields $ \{ \xi_i(t) \} $ satisfies the {\it Central Limit Theorem
(CLT) } in the Banach space $ C(T,r) $ of all continuous in the
semi-distance $ r $ functions $ f: T \to R.$ \par

  Recall that the CLT in the considered space means that for all continuous
bounded functional $ F: C(T,r) \to R $

$$
\lim_{n \to \infty} {\bf E} F(\eta_n(\cdot)) = F(\eta_{\infty}(\cdot))
$$
or equally that for all continuous functional $ F: C(T,r) \to R $

$$
\lim_{n \to \infty} Law(F(\eta_n(\cdot)) = Law( F(\eta_{\infty}(\cdot)).
$$

 Indeed, the convergence of the finite-dimensional distributions
$ \{ \eta_n(t) \} $ as $ n \to \infty $ to the finite-dimensional
distributions of a {\it Gaussian random centered continuous with probability
one relative to the distance $ r $ field } $ \eta_{\infty}(t) $
with covariation function

$$
{\bf E} \eta_{\infty}(t) \ \eta_{\infty}(s) = {\bf E} \xi_1(t) \ \xi_1(s)
$$
is evident; the {\it tightness} of the family of measures induced by the
random fields $ \{\eta_n(t) \}, \ t \in T $ in the space $ C(T,r) $ it
follows from the equality (15) for the random fields $ \eta_n(t). $ \par
 Thus, we can write, e.g., for each positive values $ u: $
$$
\lim_{n \to \infty}{\bf P} ( \sup_{t \in T} |\eta_n(t)| > u) \to
{\bf P} ( \sup_{t \in T} |\eta_{\infty}(t)| > u).
$$
 The exponential estimation (and the exact asymptotic) for the last
probability is known ([9], chapter 3).\par
 The last equality play very important role in the Monte-Carlo method and in
statistics ([9], chapter 4). \par

\vspace{3mm}

{\bf 5. Examples.} We will consider in this section a two examples random
fields where a so-called entropy integral (some generalization of Dudley's
integral, see ([6], p. 310)

$$
I = \int_0^1 \psi^{-1}( \exp H(T,d,\epsilon)) \ d \epsilon, \ \psi(x) = \exp(-\phi^*(x))
$$
 diverges. We intend to obtain in
these examples the exponential exact estimation for tail of maximum
distribution using our methods. \par
 The first example belongs to M.Talagrand [13]. \par
{\bf A. Subgaussian random field.} Let $ \{\epsilon(n) \}, n = 1,2,
\ldots, $ i.e. $ T = Z_+, $ be a sequence of independent symmetrically distributed subgaussian r.v.:

$$
{\bf P}(|\epsilon(i)| > x) = \exp(-x^2/2), \ x \ge 0,
$$
and let $ u \ge 2, $	
$$
\xi(n) = \epsilon(n)/\sqrt{ \log(n + e - 1)}.
$$
 It follows from estimation of theorem 1 after the optimization over $ C: $
$$
Q(u) \le \exp (-0.5 \ u^2 + C_0),
$$
where $ C_0 $ is some absolute constant, in the comparison to the real value
of $ Q(u), $ for which

$$
 \exp (-0.5 \ u^2) + \exp \left( - (0.5 + C_2 ) u^2 \right) \le
$$
$$
Q(u) \le \exp (-0.5 \ u^2) + \exp \left( - (0.5 + C_1 ) \ u^2 \right),
$$

(asymptotical exponential exactness). \par

\vspace{3mm}

{\bf B. Exponential bounds of distribution in the LIL for martingales.} \par

 Assume here that the martingale $ ( \xi(n), F(n) ) $ satisfies the
conditions (8a), (9) and (10). Let us choose

$$
v(n) = v_r(n) = [ \log(\log(n+3))]^{1/r},
$$
or equally

$$
v(n) = v_r(n) = [ \log(\log( \sigma(n)+3))]^{1/r},
$$
then we obtain after some calculation on the basis of theorem 1
under condition (16a) instead (16)
and choosing the {\it partition } over the "balls", more exactly, closed
intervals $ R = \{ [A(k), A(k+1) - 1] \} = \{ [A(k), B(k) ] \} $ of a view:
$$
A(k) = Q^{k-1},
$$
where $ Q = [ \ (1 + \epsilon )^k \ ] $ for $ k \ge k_0, \ \epsilon =
const > 0 $ and $ [Z] $ denotes here the integer part of $ Z; \ t_0 =
A(k), \ \delta = B(k) - A(k): $

$$
{\bf P}\left(\sup_n \frac{\xi(n)}{\sigma(n) \ v_r(n) } > u \right) \le
\exp \left[- C \ u^r \ L^{1/r}(u) \right], u > 2. \eqno(20)
$$
 In the considered case the entropy integral in general case, i.e. if

$$
\sup_n {\bf P}\left( \frac{\xi(n)}{\sigma(n) \ v_r(n)} > u \right) \ge
\exp \left[- C_0 \ u^r \ L^{1/r}(u) \right], u > 2,
$$
divergent. In detail, suppose $ \exists n_0 = 1,2,\ldots \ \Rightarrow $

$$
 {\bf P}\left( \frac{\xi(n_0)}{\sigma(n_0)} > u \right) \ge
\exp \left[- C_2 \ u^r \ L^{1/r}(u) \right], u > 2,
$$

 and let us introduce the random process (sequence)

$$
\chi(n) = \frac{ \xi(n) }{\sigma(n) \ v_r(n) },
$$
and we must add to the set $ T $ the infinite point $ \{ \infty \} $ and
define for the completeness of the set $ T: \ \chi(\infty) = 0. $ \par
We have for the natural function $ \phi_r^*(\cdot) $ for the process
$ \chi(n): $

$$
\phi_r(\lambda) \stackrel{def}{=} \log {\bf E} \sup_n \exp(\lambda \chi(n))
$$
the "tail" inequality: $ x \ge 2 \ \Rightarrow $

$$
C_1 \ x^r \ L^{1/r}(x) \le \phi_r^*(x) \le C_2 \ x^r \ L^{1/r}(x).
$$
 The natural distance $ d_{\chi} (n,m) $ for the process $ \chi(n) $
is calculated by the formula

$$
d_{\chi}(n,m) = || \chi(n) - \chi(m) ||B(\phi_r).
$$
 Put $ m = \infty; $ then we have for the amount $ N = N( T, d_{\chi}, \epsilon ) $
of optimal $ \epsilon - $ net the inequality

$$
\epsilon \ge d_{\chi}(n,\infty) \ge C / v_r(N).
$$
 We find solving the last inequality relatively the variable $ N: $

$$
H(T, d_{\chi}, \epsilon) \ge \exp( C(r) \ /\epsilon^r ), \ \epsilon \in
(0, 1/2].
$$

 The inequality [20] is in general case exact: for all the values
$ r = 2/d, \beta = d/2, \ d = 1,2,\ldots \ $
there exists a {\it polynomial} martingale $ (\xi(n), F(n)) $
satisfying the conditions (9) and (10) with $ L_1(x) = L_2(x) = L(x) = 1 $
and such that

$$
{\bf P}\left(\sup_n \frac{\xi(n)}{\sigma(n) \ v_r(n) } > u \right) \ge
\exp \left[- C_3 \ u^r \ \right], u > 2, \eqno(21a)
$$
and
$$
{\bf P} \left[\overline{\lim}_{n \to \infty} \frac{\xi(n)}{\sigma(n) \ v_r(n) } > 0 \right] > 0. \eqno(21b)
$$

 In detail, let us consider the Rademacher sequence
$ \{ \epsilon(i) \}, \ i = 1,2,\ldots; $ i.e. where
$ \{ \epsilon(i) \} $ are independent and $ {\bf P}(\epsilon(i) = 1) =
{\bf P}(\epsilon(i) = - 1) = 0.5. $ \par
  It is known that that the r. v. $ \{ \epsilon(i) \} $ belongs to the
$ B(\phi_2) $ space with the corresponding function
$$
\phi_2(\lambda) = 0.5 \ \lambda^2, \ \lambda \in (-\infty, \infty).
$$
 Indeed,
$$
 {\bf E}(\exp(\lambda \epsilon(i)) = cosh(\lambda) \le \exp(0.5 \lambda^2).
$$

 Let us denote for $ d = 1,2,3,\ldots \ \xi(n) = \xi_d(n) = $

$$
 \sum \sum \ldots \sum_{1 \le i(1)< i(2) \ldots < i(d) \le n}
 \epsilon(i(1)) \ \epsilon(i(2)) \ \epsilon(i(3)) \ldots \ \epsilon(i(d))
$$
under natural filtration $ F(n) = \sigma \{\epsilon(j), \ j \le n \}. $ \par
 It is easy to verify that $ (\xi(n), F(n)) $ is a martingale and that

 $$
        0 < C_1 \le \sigma^2(n)/n^d \le C_2 < \infty.
 $$

  We will prove the following inequality:
$$
{\bf P} \left( \overline {\lim}_{n \to \infty}
 \frac{\xi(n)}{(n \ \log ( \log(n + 3) ))^{d/2} } > 0 \right) > 0.
$$

  It is enough to consider only the case $ d = 2, $ i.e. when
 $$
 \xi(n) = \sum \sum_{1 \le i < j \le n} \epsilon(i) \ \epsilon(j).
 $$

 We observe:
 $$
 2 \ \xi(n) = \left(\sum_{k=1}^n \epsilon(k) \right)^2 - \sum_{m=1}^n (\epsilon(m))^2 \stackrel{def}{=} \Sigma_1(n) - \Sigma_2(n).
 $$

  From the classical LIL on the form belonging to Hartman-Wintner it follows
 that there exist a finite non-trivial non-negative random
variables $ \theta_1, \ \theta_2 $ for which

$$
|\Sigma_2(n)| \le n + \theta_2 \sqrt{n \ \log(\log(n + 3 )) }
$$

and

$$
\Sigma_1(n_m) \ge \theta_1 \ n_m \ \log(\log(n_m + 3))
$$

for some (random) integer positive subsequence $ n_m, \ n_m \to \infty $ as
$ m \to \infty.$ \par
 This completes the proof of inequality of (21b); the relation (21a) may be
proved by means of more fine considerations.\par

 More exactly, by means of considered method may be proved the following
relation:
$$
\overline {\lim}_{n \to \infty} \frac{\xi(n)}{(n \ \log ( \log(n + 3) ))^{d/2} }
 \stackrel{a.e}{=} \frac{2^{d/2}}{d!}.
$$

 Note that we use in the martingale case in order to estimate the variable
$ Y(u; t_0, \delta, \gamma, R) $ inside from the generic
chaining method some classical properties of martingales and $ B(\phi) $
spaces: Doob's inequality, moment estimations, connection with
$ G(\psi) \ $ norms in order to calculate the value $ Y(u; t_0,\delta, R). $\par
 Namely, let us denote $ E(k) = [A(k), B(k)].$ But we write instead the
estimation (17) for the probability $ G(u; t_0, \delta) $ the
following estimation: $ Y_k(u) \stackrel{def}{=} $

$$
Y(u; t_0,\delta, \gamma, R) \le {\bf P} \left(\max_{n \in E(k)} \xi(n) > u \
\sigma(A(k)) \ v_r(A(k))/\sigma(B(k)) \right),
$$
as long as both the functions $ \sigma(\cdot) $ and $ v_r(\cdot) $
are monotonically increasing. \par
 It follows from the Doob's inequality
$$
| \max_{n \in E(k)} \xi_n|_p \le C \ \sigma(B(k)) \cdot (p/\phi^{-1}(p))
\cdot (p/(p-1)) \le
$$
$$
2 \ C \ \sigma(B(k)) \cdot (p/\phi^{-1}(p))
$$
as long as $ p \ge 2.$ Therefore

$$
 Y_k(u) \le \exp \left(-\phi^*(C u \ \sigma(A(k))
 \ v_r(A(k))/\sigma(B(k) \right). \eqno(22)
$$
 The assertion (20) it follows from (22) after the summing over $ k. $ \par
 Moreover, if the martingale $ (\xi(n), F(n)) $ satisfies the conditions
 (8a), (9) and (10), then with probability one
$$
\overline{\lim}_{n \to \infty} \frac{\xi(n)}{\sigma(n) \ v_r(n)} \le 1,
$$
and the last inequality is exact, e.g., for the polynomial martingales [7].
\par
 Note that in this case the condition of "convergence of majoring integral"
is not satisfied.\\

 \newpage

\begin{center}

{\bf REFERENCES }\\

\end{center}

  1. {\sc Bagdasarova I.R. and Ostrovsky E.I. } (1995). A nonuniform
   exponential estimations for large deviations in a Banach space.
   {\it Theory Probab. Appl.} {\bf 45} 638-642.\\

  2. {\sc Bednorz W.} (2006). A theorem on Majorizing Measures.
    {\it Ann. Probab.} {\bf 34} 1771-1781. MR1825156 \\

  3. {\sc Fernique X.} (1975). Regularite des trajectoires des
    function aleatiores gaussiennes. {\it Ecole de Probablite de
    Saint-Flour, IV – 1974, Lecture Notes in Mathematic.} {\bf 480} 1 – 96, Springer Verlag, Berlin.\\

  4. {\sc Kozachenko Yu. V., Ostrovsky E.I.} (1985). The Banach Spaces of
      random Variables of subgaussian type. {\it Theory of Probab. and Math.
      Stat.} (in Russian). Kiev, KSU, {\bf 32}, 43 - 57.\\

  5. {\sc Kurbanmuradov O., Sabelfeld K.} (2007). Exponential bounds for
   the probability deviation of sums of random fields. Preprint.
   Weierstra$\beta$ - Institut fur Angewandte Analysis und Stochastik
   (WIAS), ISSN 0946 – 8633, p. 1-16.\\

   6. {\sc Ledoux M., Talagrand M.} (1991) Probability in Banach Spaces.
      Springer, Berlin, MR 1102015.\\

  7. {\sc Ostrovsky E.} Bide-side exponential and moment inequalities
     for tail of distribution of Polynomial Martingales. Electronic
     publication, arXiv: math.PR/0406532 v.1 Jun. 2004. \\

  8. {\sc Ostrovsky E.I.} (1999). Exponential estimations for Random Fields
     and its applications (in Russian). Russia, OINPE.\\

  9. {\sc Ostrovsky E.I.} (2002). Exact exponential estimations for random
     field maximum distribution. {\it Theory Probab. Appl.} {\bf 45} v.3,
      281 - 286. \\

   10. {\sc Talagrand M.} (1996). Majorizing measure: The generic chaining.
      {\it Ann. Probab.} {\bf 24} 1049 - 1103. MR1825156 \\

   11. {\sc Talagrand M.} (2001). Majorizing Measures without Measures.
    {\it Ann. Probab.} 29, 411-417. MR1825156 \\

   12. {\sc Talagrand M.} (2005). {\it The Generic Chaining. Upper and
     Lower Bounds of Stochastic Processes.} Springer, Berlin. MR2133757.\\

   13. {\sc Talagrand M.}(1990). Sample boundedness of stochastic processes
      under increment conditions.{\it Ann. Probab.} {\bf 18}, 1 - 49.

  \newpage

 %\vspace{5mm}

  \hspace{55mm} {\sc Department of Mathematics}\\

  \hspace{55mm} {\sc Bar-Ilan University}\\

  \hspace{55mm} {\sc Ramat Gan, Ben Gurion street, 2}\\

  \hspace{55mm} {\sc Israel \ 76521}\\

  \hspace{55mm} {\sc E-Mail}: \ galo@list.ru \\

  \hspace{55mm} {\sc E-Mail}: \ rogovee@gmail.com \\

\newpage
\begin{center}

{\bf Ostrovsky E.}\\
\vspace{4mm}

 Address: Ostrovsky E., ISRAEL, 76521, Rehovot, \ Shkolnik street.
5/8. Tel. (972)-8- 945-16-13.\\
\vspace{3mm}
e - mail: {\bf Galo@list.ru}\\

\vspace{6mm}

{\bf Rogover E.}\\
\vspace{3mm}
 Address: Rogover E., ISRAEL, 84105, Ramat Gan.

\vspace{4mm}
e - mail: {\bf rogovee@gmail.com }\\

\end{center}

\end{document}